%% file: main.tex
\documentclass[11pt,a4paper]{article}
\usepackage[margin=24mm]{geometry}
\usepackage[T1]{fontenc}
\usepackage{lmodern,microtype,amsmath,amssymb,amsthm,mathtools,booktabs,array}
\usepackage{xcolor,graphicx,caption,subcaption,tikz,rotating}
\usetikzlibrary{arrows.meta,positioning,calc,shapes.geometric,fit,backgrounds}
\usepackage[hidelinks]{hyperref}
\usepackage[numbers,sort&compress]{natbib}
\definecolor{good}{HTML}{117A75}
\definecolor{bad}{HTML}{B24A43}
\definecolor{pale}{HTML}{E9EFF2}
\definecolor{ink}{HTML}{253747}
\definecolor{blue}{HTML}{356FA1}
\newtheorem{theorem}{Theorem}
\newtheorem{lemma}{Lemma}

\newtheorem{definition}{Definition}
\newcommand{\MIS}{\mathcal M}
\newcommand{\pr}{\mathbb P}
\newcommand{\EGU}{\operatorname{EGU}}
\newcommand{\rate}{\boldsymbol\lambda}

\input{figures/visuals.tex}
\title{\vspace{-1em}\textbf{Uniformising random greedy maximal independent sets on graphs with at most five vertices}\\[5pt]\large A finite classification, a strict path obstruction, and machine-checked certification}
\author{John Fairfax-Ball\\\small Independent researcher\\\small \url{https://john-fairfaxball-research.vercel.app/}}
\date{19 November 2026}
\begin{document}
\maketitle
\begin{abstract}
A random greedy algorithm scans the vertices of a finite graph in increasing order of independent exponential arrival times and accepts each vertex unless an accepted neighbour blocks it.  We ask an inverse question: can strictly positive finite arrival rates be chosen so that the terminal \emph{inclusion-maximal} independent sets are all equally likely?  We call this exact property exponential greedy uniformisability (EGU).  We establish the complete small-order classification: every finite simple graph on at most five vertices is EGU, with the unique exception of the five-vertex path up to graph isomorphism.  A strict $2\times2$ probability minor provides an explicit obstruction for the exceptional graph; the positive cases are certified by graph-isomorphism reduction and exact constructive witnesses.  We explain the independent-exponential-clock probability law and the formal Lean 4 verification, registered as Palomar \texttt{PALOMAR-2026-11-19-000021} (version 1).  The registration is a formal proof record, not a claim of novelty certification or journal peer review.
\end{abstract}

\section{Introduction}
Random sequential adsorption on graphs (RSA) and random greedy independent-set algorithms start from a random ordering and accept a vertex when no previously accepted neighbour forbids it \citep{Pippenger1989,Krivelevich2024}.  The result is maximal by inclusion, but generally neither maximum by cardinality nor uniformly distributed among maximal independent sets.  An equal-rate exponential race gives a uniform random permutation.  Here, in contrast, the rates are design parameters: we wish to \emph{control the entire terminal-set law}.  Unequal rates in RSA are not themselves new \citep{PenroseRosoman2012}; the question is whether they can enforce an \emph{exact uniform law} over the finite set of terminal configurations.

\begin{figure}[htbp]
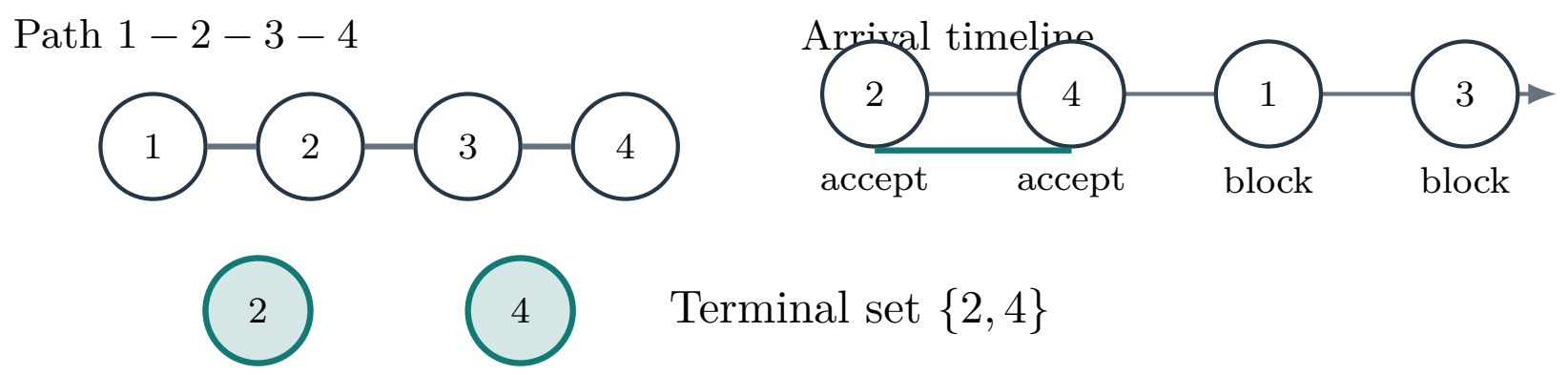
\centering\figscan
\caption{\textbf{One deterministic outcome of the exponential-clock scan.}  On the path $1-2-3-4$, the illustrative strict arrival order is $2,4,1,3$ (not a claim about observed frequencies).  Vertices 2 and 4 are accepted; 1 and 3 are blocked by accepted neighbours.  The terminal set $\{2,4\}$ is inclusion-maximal.  Filled teal circles are accepted; crossed grey circles are rejected.}\label{fig:scan}\end{figure}

\begin{figure}[htbp]
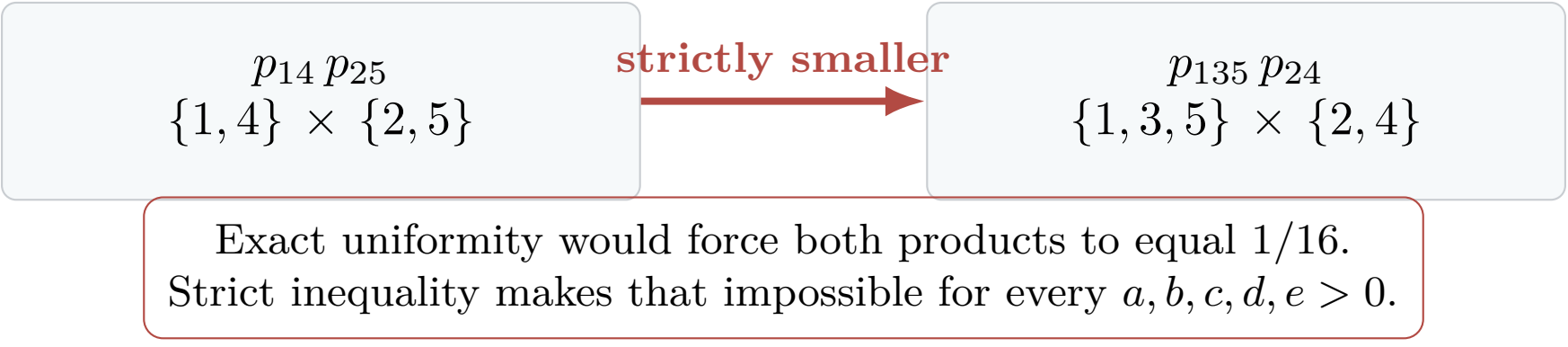
\centering\figobstruction
\caption{\textbf{The $P_5$ obstruction is algebraic, not numerical.}  For any positive rates, the product of the output masses on the lower pair of sets strictly exceeds the product on the upper pair.  Equal mass $1/4$ for all four would require both products to equal $1/16$, contradicting the strict inequality.  Equation~\eqref{eq:delta} supplies the positive polynomial certificate.}\label{fig:obstruction}\end{figure}

\begin{figure}[htbp]
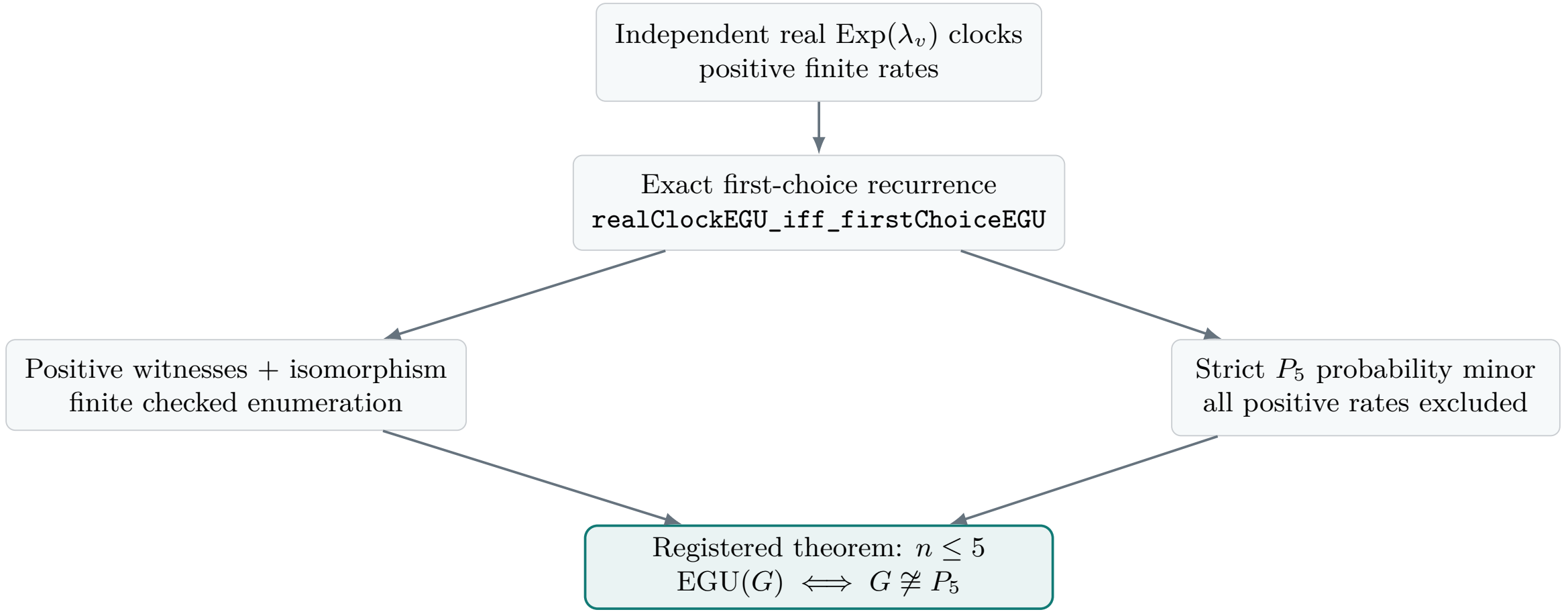
\centering\scalebox{.87}{\figarchitecture}
\caption{\textbf{Logical dependency structure (conceptual, not a literal import graph).}  The genuine independent exponential-clock measure is linked by a proved law to the exact finite first-choice recursion.  Positive representative certificates plus orbit/isomorphism reduction establish the constructive cases; the strict $P_5$ minor excludes the exception.  The branches meet at the registered $n\le5$ biconditional.}\label{fig:architecture}\end{figure}

This distinction is important.  Equal-rate output uniformity on trees was studied in \citet{FairfaxBall2026}; allowing arbitrary positive rates changes the answer even for the three-vertex path (Figure~\ref{fig:p3}).  Uniform maximal-set sampling is also different from Gibbs distributions and maximum-weight independent-set optimisation.  We give a small, complete classification, not a classification of arbitrary finite graphs.  Figures~\ref{fig:scan}--\ref{fig:architecture} are explanatory vector diagrams whose exact claims have separate source and validation records.

\section{Model and the first-choice identity}
Let $G=(V,E)$ be a finite simple graph.  Independently for each $v\in V$, sample $T_v\sim\mathrm{Exp}(\lambda_v)$, where $0<\lambda_v<\infty$.  Inspect vertices in increasing $T_v$ order, accepting $v$ precisely if it has no previously accepted neighbour.  The output $I$ is independent.  A vertex omitted from $I$ was blocked when inspected, so it has a neighbour in $I$: thus $I$ is \emph{inclusion-maximal}.  Write $\MIS(G)$ for the family of all such sets, and $p_{G,\rate}(I)$ for the output probability.

\begin{definition}\label{def:egu}
A graph $G$ is \emph{exponentially greedily uniformisable} (EGU) if some vector $\rate\in(0,\infty)^V$ satisfies
\[p_{G,\rate}(I)=\frac{1}{|\MIS(G)|}\qquad\text{for every }I\in\MIS(G).\]
When $V=\varnothing$, the unique output is the empty set, and $G$ is EGU.
\end{definition}

The exponential-race identity says that the first clock belongs to $v$ with probability $\lambda_v/\Lambda$, where $\Lambda=\sum_u\lambda_u$.  If $v$ is accepted first, $v$ and its neighbours are removed from subsequent consideration.  Memorylessness yields the following exact recursion (equivalently, sum over permutations weighted by exponential-race probabilities):
\begin{equation}\label{eq:race}
p_{G,\rate}(I)=\sum_{v\in I}\frac{\lambda_v}{\Lambda}\;p_{G-N[v],\rate|_{V\setminus N[v]}}(I\setminus\{v\}),
\end{equation}
where $N[v]$ denotes the closed neighbourhood, $p_{\varnothing}(\varnothing)=1$, and the probability is defined as zero for an impossible target.  The identity uses the \emph{actual independent exponential clocks}, not a replacement sampling model.  The formal bridge appears in \texttt{RealClockGreedyLaw.lean}.

\section{Positive-rate control: two small examples}
\begin{figure}[htbp]
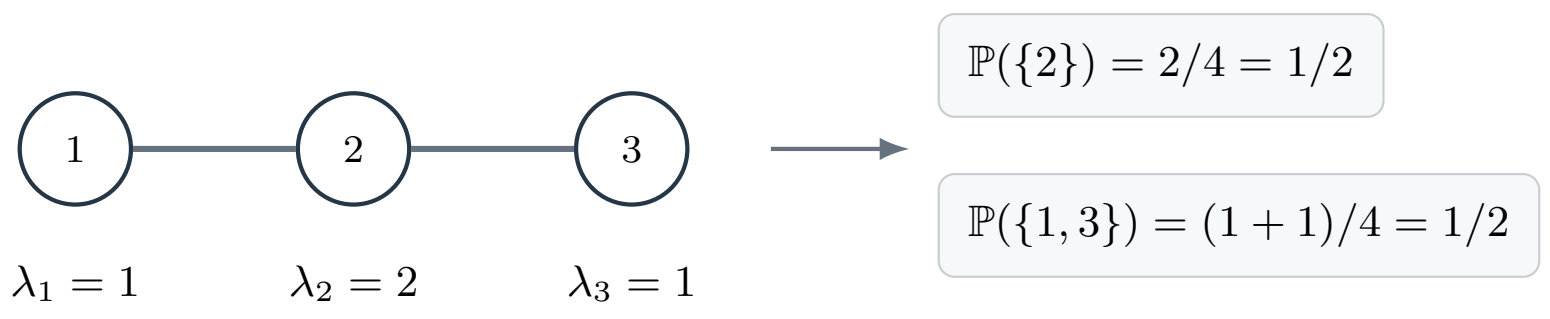
\centering\figpthree
\caption{\textbf{Rate design on $P_3$.}  The vertex labels are above the path and rates below it.  With rates $(1,2,1)$, the centre clock wins with probability $2/4$, giving $\{2\}$; either endpoint wins with total probability $2/4$, giving $\{1,3\}$.  Both are exact fractions, not numerical estimates.}\label{fig:p3}\end{figure}

For the path $P_3=1-2-3$, the only maximal independent sets are $\{2\}$ and $\{1,3\}$.  Whichever vertex arrives first determines which of these two outputs is produced:
\[p(\{2\})=\frac{\lambda_2}{\lambda_1+\lambda_2+\lambda_3},\qquad p(\{1,3\})=\frac{\lambda_1+\lambda_3}{\lambda_1+\lambda_2+\lambda_3}.\]
The strictly positive rates $(1,2,1)$ therefore give $1/2$ to each outcome (Figure~\ref{fig:p3}).  Choosing all rates equal would instead give $1/3$ and $2/3$.

The path $P_4=1-2-3-4$ is instructive because it is not a cograph, yet is EGU.  Take reflection-symmetric rates
\[(\lambda_1,\lambda_2,\lambda_3,\lambda_4)=(1,\sqrt3-1,\sqrt3-1,1).\]
By applying \eqref{eq:race},
\[p(\{1,4\})=\frac{1}{(1+\sqrt3-1)^2}=\frac13,\qquad p(\{1,3\})=p(\{2,4\})=\frac13.\]
These are positive finite real rates (not a limiting construction).  The graph and corresponding complementary witness types are visualised in Figure~\ref{fig:witnesses}.  A compact algebraic verification is in Appendix~\ref{app:witnesses}.

\section{The complete small-order classification}
\begin{theorem}[At-most-five-vertex classification]\label{thm:main}
For every finite simple graph $G$ with $|V(G)|\le5$,
\[\EGU(G)\quad\Longleftrightarrow\quad G\not\cong P_5.\]
In particular, the exception is \emph{whole-graph isomorphism} to the five-vertex path, not merely the appearance of an induced $P_5$ in a larger graph.
\end{theorem}

\begin{figure}[htbp]
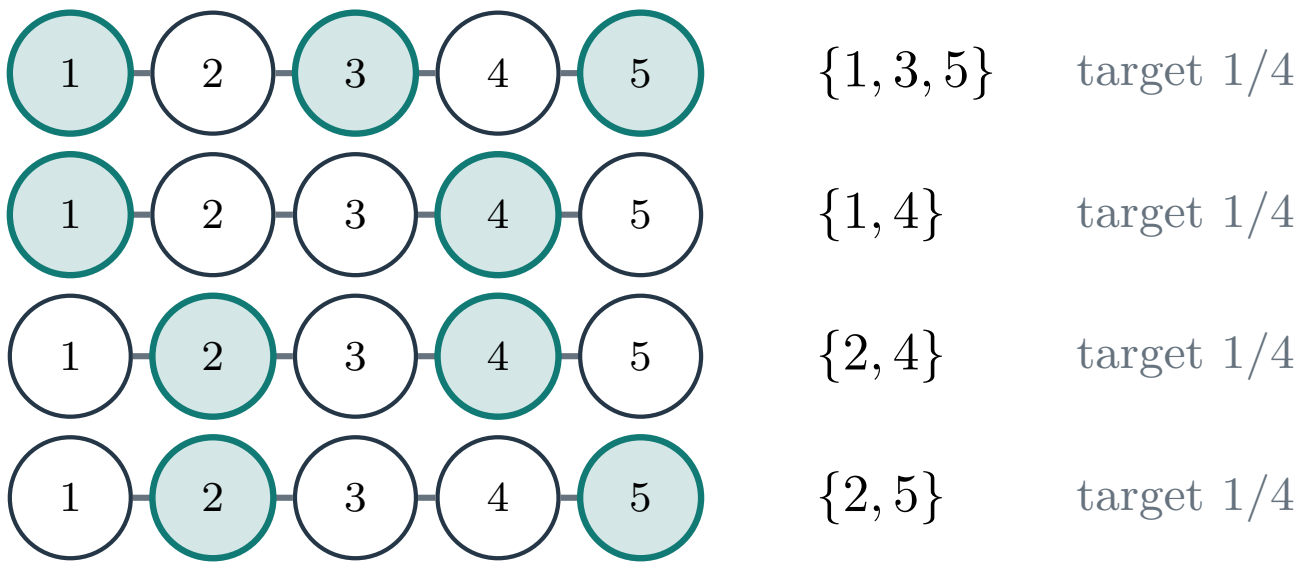
\centering\figpfivesets
\caption{\textbf{All four terminal sets of $P_5$.}  The four rows show $\{1,3,5\}$, $\{1,4\}$, $\{2,4\}$ and $\{2,5\}$; teal filled vertices belong to the indicated set.  The exact-uniform target is $1/4$ for each row.  These four sets were checked by exhaustive independent-set enumeration.}\label{fig:p5sets}\end{figure}

The $P_5$ graph has four maximal independent sets, drawn in Figure~\ref{fig:p5sets}.  The required target law is thus $1/4$ for each of them; one cannot merely balance the probability of output cardinalities.

\subsection{A strict obstruction for $P_5$}
Label the path $1-2-3-4-5$, and abbreviate the positive rates by $(a,b,c,d,e)$; let $S=a+b+c+d+e$.  Define $p_{135}=p(\{1,3,5\})$, and similarly for the other three outputs.  Applying the first-choice identity gives
\begin{align}\label{eq:p5prob}
p_{14}&=\frac{ad}{(a+b)(c+d+e)}, &p_{25}&=\frac{be}{(a+b+c)(d+e)},\\
p_{24}&=\frac{bd}{S}\left(\frac1{a+b}+\frac1{d+e}\right), &
p_{135}&=\frac1S\left(\frac{a(c+e)}{c+d+e}+c+\frac{e(a+c)}{a+b+c}\right).\notag
\end{align}

\begin{lemma}[Strict output minor]\label{lem:strict}
For all $a,b,c,d,e>0$, the $P_5$ output law obeys
\begin{equation}\label{eq:minor}p_{14}p_{25}<p_{135}p_{24}.\end{equation}
\end{lemma}
\begin{proof}
Put $A=a+b$, $B=d+e$, $L=a+b+c$, and $R=c+d+e$.  Direct substitution into \eqref{eq:p5prob} and clearing only strictly positive denominators gives
\begin{equation}\label{eq:delta}
p_{135}p_{24}-p_{14}p_{25}=\frac{bcd\,Q(a,b,c,d,e)}{ABLR\,S^2},
\end{equation}
where the explicit polynomial $Q$ is displayed as a sum of polynomials with \emph{strictly positive coefficients} in Appendix~\ref{app:polynomial}.  Thus $Q>0$, and every other factor on the right side of \eqref{eq:delta} is positive.  Equality is impossible for finite strictly positive rates.
\end{proof}

The strict minor proves $P_5$ cannot be EGU.  This argument rules out \emph{all} real positive rate choices; a failed search over rational candidates, a numerical minimum or a limiting sequence would not establish the same conclusion.

\subsection{The other isomorphism classes}
\begin{sidewaysfigure}[p]\centering\includegraphics[width=.88\textheight]{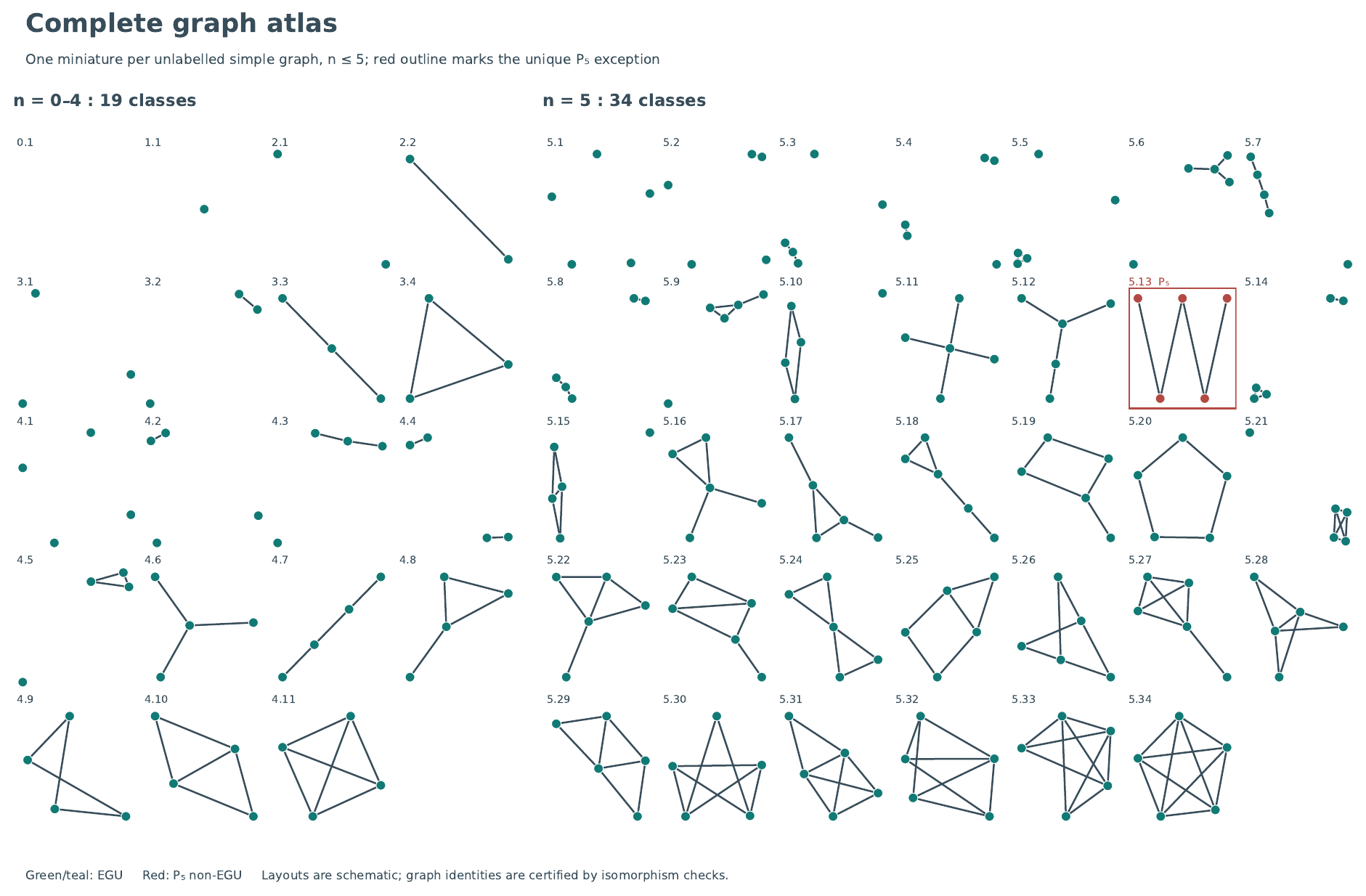}
\caption{\textbf{The complete unlabelled graph atlas for $0\le n\le5$.}  Each miniature is a distinct isomorphism class (53 in total).  The one red-highlighted class is $P_5$; the other 52 are EGU.  Positions are schematic, not a geometric embedding claim.  Generated deterministically by \texttt{figures/generate\_atlas.py} from NetworkX's \texttt{graph\_atlas\_g}, with a separate independent isomorphism/count check; the formal completeness proof uses its own Lean edge-vector enumeration.}\label{fig:atlas}\end{sidewaysfigure}

For the positive direction, the formal development represents each graph by its $\binom n2$-bit edge vector, canonically compares its isomorphism class to a fixed finite representative table, and checks an exact rate witness for each nonexceptional representative.  The finite case split is kernel-checked by \texttt{allBool\_sound}, and the transfer to arbitrary labellings uses \path{realClockEGU_congr_equiv}.  The computed classes are summarised in Table~\ref{tab:classes}; a fully drawn atlas is Figure~\ref{fig:atlas}.  The tables are part of the formal proof, not an unexplained numerical optimisation.
\begin{table}[htbp]\centering\caption{Unlabelled simple-graph classes.  All positive classes possess exact positive witnesses; $P_5$ is the sole negative class within the registered scope.}\label{tab:classes}\begin{tabular}{rrrr}\toprule $n$ & Isomorphism classes & EGU & Non-EGU\\\midrule 0&1&1&0\\1&1&1&0\\2&2&2&0\\3&4&4&0\\4&11&11&0\\5&34&33&1\\\midrule Total&53&52&1\\\bottomrule\end{tabular}\end{table}

\begin{figure}[htbp]
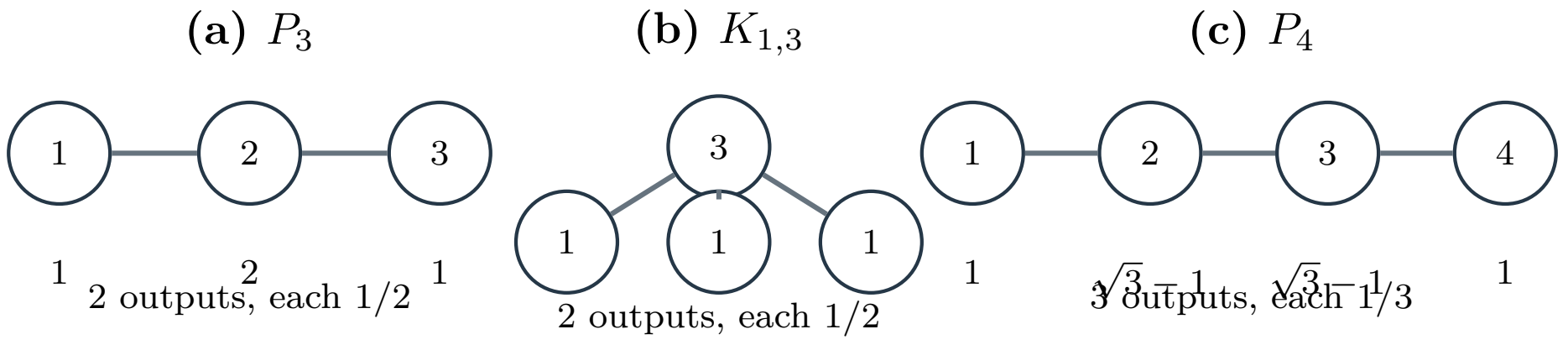
\centering\figwitnesses
\caption{\textbf{Three exact constructive witnesses.}  (a) $P_3$ with rates $(1,2,1)$ has two equal output masses; (b) $K_{1,3}$ with hub rate 3 and leaves of rate 1 gives mass $1/2$ to the hub and $1/2$ to the set of all leaves; (c) $P_4$ with rates $(1,\sqrt3-1,\sqrt3-1,1)$ gives each of its three maximal independent sets mass $1/3$.  Numbers by vertices are rates, not probabilities.}\label{fig:witnesses}\end{figure}

\begin{proof}[Proof of Theorem~\ref{thm:main} (certified finite classification)]
For $n\le4$, all $1+1+2+4+11$ unlabelled classes are covered by exact positive-rate witnesses and the checked isomorphism reduction.  For $n=5$, the same certified enumeration finds 34 unlabelled classes: 33 carry witnesses, and the remaining class is isomorphic to $P_5$.  Those certificates establish the positive implication in every labelled graph through five vertices.  Lemma~\ref{lem:strict} rules out $P_5$, and graph-isomorphism invariance transports that obstruction to every labelling of this path.  The empty graph is handled separately.  The full finite case reduction, witness verification, and genuine-clock link are checked in Lean (see Section~\ref{sec:formal}).
\end{proof}

\section{Formal verification and evidence tiers}\label{sec:formal}
The public Challenge states the biconditional of Theorem~\ref{thm:main} for all symmetric irreflexive neighbour maps on $\mathrm{Fin}(n)$, $n\le5$, with a positive rate vector and actual independent $\mathrm{Exp}(\lambda_v)$ product measure.  The imported Solution combines the real-clock law \path{realClockEGU_iff_firstChoiceEGU}, the finite classification \path{graphs_le_five_classification}, the strict obstruction \path{p5Nb_not_realClockEGU}, and isomorphism transfer.  The Challenge's intentional placeholder is not a proof hole in the accepted Solution.

At the pinned formal submission, the comparator and independent Lean, NanoDa and con-ron checks accepted the solution; the recorded axiom dependency list is \texttt{propext}, \texttt{Classical.choice}, and \texttt{Quot.sound}.  The public identifier supplied by the owner is \href{https://palomar-registry.org/entry.html?id=PALOMAR-2026-11-19-000021\&version=1}{PALOMAR-2026-11-19-000021 (v1)}.  The cited public entry records the formal verification of the statement and proof at a pinned source snapshot.  This paper does not treat registration as independent editorial refereeing.

\section{Related work, scope and open directions}
Pippenger's RSA analysis introduced an influential graph-level sequential occupancy model \citep{Pippenger1989}.  Krivelevich et al. study random greedy maximal independent sets and local-limit densities, with primarily uniform orderings \citep{Krivelevich2024}.  Rate-asymmetric RSA is known \citep{PenroseRosoman2012}, so neither exponential arrivals nor unequal rates should be presented as our discovery.  The equal-rate tree-uniformity results in \citet{FairfaxBall2026} address a different quantifier.  The precise distinctive issue studied here is whether \emph{some} positive fixed per-vertex rate vector makes the \emph{full output law} exactly uniform.

The theorem deliberately stops at five vertices.  It does not show that all graphs containing an induced $P_5$ fail EGU, classify $n=6$, or assert that approximate uniformisation implies exact uniformisation.  Work in the repository also studies cograph closure, unions and joins, and a six-vertex census; those are separate claims and are not promoted to parts of the registered small-order iff.  Future work includes structural characterisations beyond five vertices, quantitative distance to uniformity for forbidden examples, and rate-design algorithms with rigorous feasibility certificates.

\paragraph{Methods and disclosure.} Research direction and authorship are attributed to John Fairfax-Ball.  The development and this draft used AI-assisted research, Lean code, scripts and editorial assistance.  Exact certificates and independently checked kernel replay supply the mathematical assurance claimed here; no human referee review or comprehensive novelty certification is implied.  The formal repository and figure tests are publicly auditable.

\appendix
\sloppy
\section{The positive polynomial certificate}\label{app:polynomial}
For compactness, group $Q$ in powers of the middle rate $c$:
\[Q=Q_0+cQ_1+c^2(a+b+d+e),\]
with
\begin{align*}
Q_1={}&2a^2+3ab+3ad+3ae+b^2+2bd+3be+d^2+3de+2e^2,\\
Q_0={}&a^3+2a^2(b+d+e)\\
&+a(b^2+3bd+3be+d^2+3de+2e^2)\\
&+b^2(d+e)+b(d^2+3de+2e^2)\\
&+d^2e+2de^2+e^3.
\end{align*}
Every monomial coefficient in $Q_0,Q_1$ and $a+b+d+e$ is strictly positive, so $Q>0$ for strictly positive $a,b,c,d,e$.  The equality \eqref{eq:delta} can be checked by rational-function simplification of \eqref{eq:p5prob}; \texttt{figures/validate.py} independently recomputes the factorisation with SymPy.  The same sign argument is certified in Lean in \texttt{SmallGraphObstructions.lean}, under \texttt{p5Nb\_strict\_minor}.

\section{A direct $P_4$ verification}\label{app:witnesses}
Write the symmetric rates as $(a,b,b,a)$ with $a,b>0$.  The first-choice recursion shows that the probability of the endpoint pair $\{1,4\}$ is $a^2/(a+b)^2$.  By reflection symmetry the remaining two maximal independent sets have equal mass, and all three masses sum to 1.  Taking $a/(a+b)=1/\sqrt3$, i.e. $(a,b)=(1,\sqrt3-1)$, makes all masses $1/3$.  The witness therefore uses genuine positive rates, not an approximation.

\section{Reproduction and certificate trace}\label{app:reproduce}
The formal code and independent figure-check scripts are publicly available at \url{https://github.com/jfairfaxball-348/Inverse-Greedy}. The repository figure ledger records sources and validation. The core formal interface is in \path{Challenge.lean} and \path{Solution.lean}; classification, obstruction and real-clock law are imported from the \path{InverseGreedy/} directory. The independent figure checker is not a substitute for the registered formal verification.

\bibliographystyle{plainnat}
\bibliography{references}
\end{document}

%% file: figures/visuals.tex
\tikzset{
 igvertex/.style={circle,draw=ink,thick,minimum size=8mm,inner sep=0pt,fill=white,font=\scriptsize},
 igselected/.style={igvertex,fill=good!18,draw=good,very thick},
 igblocked/.style={igvertex,fill=pale,draw=ink!45},
 igedge/.style={line width=1.15pt,draw=ink!70},
 igbox/.style={draw=ink!25,rounded corners=3pt,fill=pale!40,inner sep=6pt},
 igarrow/.style={-{Latex[length=2mm]},thick,draw=ink!70}
}
\newcommand{\figscan}{%
\begin{tikzpicture}[x=1cm,y=1cm,scale=.91,every node/.style={transform shape}]
\node[anchor=west,font=\footnotesize] at (0,1.0) {Path $1-2-3-4$};
\foreach \i in {1,...,4} {\node[igvertex] (g\i) at (1.2*\i,.15) {\i};}
\foreach \i/\j in {1/2,2/3,3/4} \draw[igedge] (g\i)--(g\j);
\node[anchor=west,font=\footnotesize] at (6.0,1.0) {Arrival timeline};
\draw[igarrow] (6.4,.55)--(11.9,.55);
\foreach \v/\x/\lab in {2/6.7/accept,4/8.2/accept,1/9.7/block,3/11.2/block} {
  \node[igvertex,fill=white] at (\x,.55) {\v};
  \node[font=\scriptsize,anchor=north] at (\x,.13) {\lab};}
\draw[good,very thick] (6.7,.12)--(8.2,.12);
\node[igselected] at (2.0,-1.1) {2};\node[igselected] at (4.0,-1.1) {4};
\node[font=\small,anchor=west] at (5.0,-1.1) {Terminal set $\{2,4\}$};
\end{tikzpicture}}
\newcommand{\figpthree}{%
\begin{tikzpicture}[scale=.92,every node/.style={transform shape}]
\foreach \n/\x/\r in {1/0/1,2/2/2,3/4/1} {
 \node[igvertex] (p\n) at (\x,0) {\n};
 \node[font=\footnotesize,below=7mm] at (p\n) {$\lambda_{\n}=\r$};}
\draw[igedge] (p1)--(p2)--(p3);
\draw[igarrow] (5.0,0)--(6.0,0);
\node[igbox,anchor=west,align=left,font=\footnotesize] at (6.2,.6) {$\pr(\{2\})=2/4=1/2$};
\node[igbox,anchor=west,align=left,font=\footnotesize] at (6.2,-.55) {$\pr(\{1,3\})=(1+1)/4=1/2$};
\end{tikzpicture}}
\newcommand{\figpfivesets}{%
\begin{tikzpicture}[x=.82cm,y=.78cm]
\foreach \y/\vals/\lab in {0/{1,3,5}/{$\{1,3,5\}$},-1.2/{1,4}/{$\{1,4\}$},-2.4/{2,4}/{$\{2,4\}$},-3.6/{2,5}/{$\{2,5\}$}} {
 \foreach \i in {1,...,5}{
  \node[igvertex] (n\i) at (\i*1.15,\y) {\i};
 }
 \foreach \a/\b in {1/2,2/3,3/4,4/5}{\draw[igedge] (n\a)--(n\b);}
 \foreach \v in \vals {\node[igselected] at (\v*1.15,\y) {\v};}
 \node[font=\small,anchor=west] at (7.0,\y) {\lab};
 \node[font=\footnotesize,anchor=west,color=ink!70] at (9.1,\y) {target $1/4$};
 }
\end{tikzpicture}}
\newcommand{\figobstruction}{%
\begin{tikzpicture}[font=\small,node distance=1cm]
\node[igbox,align=center,minimum width=4.7cm,minimum height=1.45cm] (lhs) {$p_{14}\,p_{25}$\\ $\{1,4\}\,\times\,\{2,5\}$};
\node[igbox,align=center,minimum width=4.7cm,minimum height=1.45cm,right=2.1cm of lhs] (rhs) {$p_{135}\,p_{24}$\\ $\{1,3,5\}\,\times\,\{2,4\}$};
\draw[-{Latex[length=3mm]},line width=1.5pt,color=bad] (lhs)-- node[above,font=\bfseries\footnotesize,text=bad] {strictly smaller} (rhs);
\node[draw=bad,rounded corners,align=center,below=7mm of $(lhs)!0.5!(rhs)$,inner sep=5pt,font=\footnotesize] {Exact uniformity would force both products to equal $1/16$.\\ Strict inequality makes that impossible for every $a,b,c,d,e>0$.};
\end{tikzpicture}}
\newcommand{\figwitnesses}{%
\begin{tikzpicture}[x=1cm,y=1cm]
\node[font=\small\bfseries] at (1.5,1.35) {(a) $P_3$};
\foreach \x/\l/\r in {0/1/1,1.5/2/2,3/3/1}{\node[igvertex] (a\l) at (\x,.4) {\l};\node[below=3mm of a\l,font=\scriptsize] {$\r$};}
\draw[igedge] (a1)--(a2)--(a3);
\node[font=\scriptsize] at (1.5,-.75) {2 outputs, each $1/2$};
\node[font=\small\bfseries] at (5.2,1.35) {(b) $K_{1,3}$};
\node[igvertex] (bh) at (5.2,.45) {$3$};
\foreach \x/\idx in {4.0/1,5.2/2,6.4/3}{\node[igvertex] (b\idx) at (\x,-.3) {$1$};\draw[igedge] (bh)--(b\idx);}
\node[font=\scriptsize] at (5.2,-.9) {2 outputs, each $1/2$};
\node[font=\small\bfseries] at (9.4,1.35) {(c) $P_4$};
\foreach \x/\i/\r in {7.2/1/1,8.6/2/{\sqrt3-1},10/3/{\sqrt3-1},11.4/4/1}{\node[igvertex] (c\i) at (\x,.4) {\i};\node[below=3mm of c\i,font=\scriptsize] {$\r$};}
\draw[igedge] (c1)--(c2)--(c3)--(c4);
\node[font=\scriptsize] at (9.4,-.75) {3 outputs, each $1/3$};
\end{tikzpicture}}
\newcommand{\figarchitecture}{%
\begin{tikzpicture}[node distance=6mm and 6mm,font=\footnotesize]
\node[igbox,align=center,minimum width=5cm] (clock) {Independent real $\mathrm{Exp}(\lambda_v)$ clocks\\positive finite rates};
\node[igbox,align=center,below=of clock,minimum width=5cm] (law) {Exact first-choice recurrence\\\texttt{realClockEGU\_iff\_firstChoiceEGU}};
\node[igbox,align=center,below left=10mm and 12mm of law,minimum width=4.4cm] (pos) {Positive witnesses + isomorphism\\finite checked enumeration};
\node[igbox,align=center,below right=10mm and 12mm of law,minimum width=4.4cm] (neg) {Strict $P_5$ probability minor\\all positive rates excluded};
\node[draw=good,fill=good!9,rounded corners,thick,align=center,below=17mm of law,minimum width=5.3cm,yshift=-1.4cm] (out) {Registered theorem: $n\le5$\\ $\EGU(G)\iff G\not\cong P_5$};
\draw[igarrow] (clock)--(law);
\draw[igarrow] (law)--(pos);\draw[igarrow] (law)--(neg);
\draw[igarrow] (pos)--(out);\draw[igarrow] (neg)--(out);
\end{tikzpicture}}

%% file: main.bbl
\begin{thebibliography}{4}
\providecommand{\natexlab}[1]{#1}

\bibitem[Fairfax-Ball(2026)]{FairfaxBall2026}
John Fairfax-Ball.
\newblock Greedy uniformity on trees: Exact obstruction and near-uniform spiders.
\newblock \emph{arXiv preprint arXiv:2610.02276}, 2026.
\newblock \url{https://arxiv.org/abs/2610.02276}.

\bibitem[Krivelevich et~al.(2024)Krivelevich, M{\'e}sz{\'a}ros, Michaeli, and Shikhelman]{Krivelevich2024}
Michael Krivelevich, Tam{\'a}s M{\'e}sz{\'a}ros, Peleg Michaeli, and Clara Shikhelman.
\newblock Greedy maximal independent sets via local limits.
\newblock \emph{Random Structures \& Algorithms}, 64(4):986--1015, 2024.
\newblock doi:10.1002/rsa.21200.

\bibitem[Penrose and Rosoman(2012)]{PenroseRosoman2012}
Mathew D. Penrose and Thomas Rosoman.
\newblock Percolation of even sites for random sequential adsorption.
\newblock \emph{Stochastic Processes and their Applications}, 122(4):1866--1886, 2012.
\newblock doi:10.1016/j.spa.2012.01.007.

\bibitem[Pippenger(1989)]{Pippenger1989}
Nicholas Pippenger.
\newblock Random sequential adsorption on graphs.
\newblock \emph{SIAM Journal on Discrete Mathematics}, 2(3):393--401, 1989.
\newblock doi:10.1137/0402034.

\end{thebibliography}
